\documentclass[times,doublespace]{article}
\usepackage{latexsym,amssymb,amsfonts}
\usepackage{amsmath}

\newtheorem{theorem}{Theorem}
\newtheorem{corollary}{Corollary}

\begin{document}
\title{All admissible meromorphic solutions of Hayman's equation}
%\shorttitle{Meromorphic solutions of Hayman's equation}

% Enter the publication year and the ID number of the paper
%\volumeyear{XXX}
%\paperID{XXX}

% Author name(s)
\author{Rod Halburd\footnote{Department of Mathematics, University College London, 
Gower Street,
London WC1E 6BT, UK}\ \ and Jun Wang\footnote{School of Mathematical Sciences, Fudan University, 
Shanghai 200433, China}}
% Abbreviated author name for running headers
%\abbrevauthor{R. Halburd and J. Wang}
% Abbreviated author name for first page header
%\headabbrevauthor{R. Halburd and J. Wang}

% Address / e-mail address of corresponding author
%\correspdetails{R.Halburd@ucl.ac.uk}

% Enter received/revised/accepted dates as necessary
%\received{XXX}
%\revised{XXX}
%\accepted{XXX}

% Enter details of editor communicating this article
%\communicated{A. Editor}
%\author{Rod Halburd and Jun Wang}

\maketitle

\begin{abstract}
We find all non-rational meromorphic solutions of the equation
$ww''-(w')^2=\alpha(z)w+\beta(z)w'+\gamma(z)$,
where $\alpha$, $\beta$ and $\gamma$ are rational functions of $z$.
In so doing we answer a question of Hayman by showing that all such
solutions have finite order. Apart from special choices of the
coefficient functions, the general solution is not meromorphic and
contains movable branch points.  For some choices for the
coefficient functions the equation admits a one-parameter family of
non-rational meromorphic solutions.  Nevanlinna theory is used to
show that all such solutions have been found and allows us to avoid
issues that can arise from the fact that resonances can occur at
arbitrarily high orders.  We actually solve the more general problem of
finding all  meromorphic solutions that are admissible in the sense of Nevanlinna theory, where the coefficients
$\alpha$, $\beta$ and $\gamma$ are meromorphic functions.

%\vskip 5mm

%\noindent{\it{Keywords:}} Complex differential equations, Meromorphic
%solutions, Nevanlinna theory.

%\noindent {\it 2010 Mathematics Subject Classification:}
%34M10, 30D35.

\end{abstract}

\section{Introduction}
Local series methods often provide strong necessary conditions for
the general solution of an ordinary differential equation to have a
meromorphic general solution.  The existence of a meromorphic
general solution (or more generally, that an ODE has the Painlev\'e
property, see e.g.~\cite{ablowitzc:91}) is often used as a way to identify
equations that are integrable, i.e., in some sense exactly solvable.
We wish to extend this idea to that of finding all sufficiently
complicated meromorphic solutions of an ODE, even when the general
solution is not meromorphic.  We are effectively using singularity
structure to find integrable sectors of the solution space of
the equation under consideration.

In this paper we
will find all {\em admissible} meromorphic solutions of the
differential equation
\begin{equation}
\label{he}
ww''-w'^2=\alpha(z)w+\beta(z)w'+\gamma(z) ,
\end{equation}
where $\alpha$, $\beta$ and $\gamma$ are meromorphic functions.
Heuristially a meromorphic solution is admissible if it is more
complicated than the coefficients that appear in the equation.  In
particular, if the coefficients are rational functions then any
transcendental (i.e. non-rational) meromorphic solution is
admissible.  If the coefficients are constants then any non-constant
meromorphic solution is admissible. The precise definition of an
admissible meromorphic solution $w$ of Eq.~\eqref{he} is that $w$
satisfies
\begin{equation}
\label{admissible}
T(r,\alpha)+T(r,\beta)+T(r,\gamma)=S(r,w),
\end{equation}
where $T$ is the Nevanlinna characteristic and $S(r,w)$ is used to
denote any function of $r$ that is $o(T(r,w))$ as $r\to\infty$
outside of some possible exceptional set of finite linear measure.

The main result of this paper is the following.
%\vskip 2mm \noindent {\bf Theorem 1.}\,\,{\em
\begin{theorem}
\label{mainthm}
Suppose that $w$ is a meromorphic solution of Eq.~\eqref{he}, where
the meromorphic coefficients $\alpha(z)$, $\beta(z)$ and $\gamma(z)$
satisfy Eq.~\eqref{admissible}.
% If the common zeros of $\alpha(z),\beta(z)$ and $\gamma(z)$ are finitely many, t
Then $w$ is one of the solutions described in the following list, where $c_1$ and $c_2$ are constants.
\begin{enumerate}
%\item  If $\alpha\equiv\beta\equiv\gamma\equiv 0$, then $w(z)=c_2{\rm e}^{c_1z}$.
\item \label{concA}If $\beta\equiv\gamma\equiv 0$ and $k_1=\alpha\not\equiv 0$ is a constant, then
$w=\frac{k_1}{c_1^2}\left\{1+\cosh(c_1z+c_2)\right\}$ or $w=-\frac{k_1}2(z+c_2)^2$.
\item \label{concB} If $\gamma\equiv 0$, $\beta\not\equiv 0$ and $k_1=-\alpha/\beta$ is a constant, then $w(z)=c_1{\rm e}^{k_1z}$.
\item \label{concC} If $\gamma\equiv 0$ and $\alpha+\beta'\equiv 0$, then $w(z)={\rm e}^{c_1z}\left\{c_2-\displaystyle\int\beta(z){\rm e}^{-c_1z}{\rm d}z\right\}$.
\item \label{concD} If $\gamma\not\equiv 0$ and there is a constant $k_1$ and a meromorphic function $h$ satisfying
$h^2+\beta h+\gamma=0$
and
$h'-k_1h=\alpha+k_1\beta$, then
$w={\rm e}^{k_1 z}\left(c_1+\int h(z){\rm e}^{-k_1 z}{\rm d}z\right)$.
\item \label{concE}
Suppose that $\gamma\not\equiv 0$ and $A=\frac{\beta(\alpha+\beta')-\gamma'}{\gamma}$ is a constant. 
% Let
%$$
%h(z)=\left(
%\frac\beta 2A-\beta'-2\alpha
%\right){\rm e}^{Az/2}.
%$$
\begin{enumerate}
\item If $A=0$ and there is a nonzero constant $k_1$ such that
$k_1^2\beta+\beta''+2\alpha'=0$, then
$$
k_2^2
=
\frac1{k_1^2}
\left\{
\frac{1}{4k_1^2}\left(
\beta'+2\alpha
\right)^2+\left(\gamma-\frac{\beta^2}4\right)
\right\}
$$
is also a constant.  If $k_2\ne0$
then
$w=\pm k_2\cosh(k_1 z+c_1)+\frac {\beta'+2\alpha}{2k_1^2}$.
\item
If $k_1^2=\frac{\left(\frac\beta 2 A-\beta'-2\alpha\right)^2}{\beta^2-4\gamma}$ is a nonzero constant then 
$w=c_1{\rm e}^{\left(-\frac A2\pm k_1\right)z}-\frac 1{2k_1^2}\left(\frac\beta 2A-\beta'-2\alpha\right)$.
%{\rm e}^{-Az/2}\left\{
%c_1{\rm e}^{\pm k_1z}-\frac h{2k_1^2}
%\right\}$ 
\item
If $\alpha$ and $\gamma$ are non-zero constants and $\beta=0$, then $w(z)=-\frac\alpha 2(z+c_1)^2-\frac{\gamma}{2\alpha}$.
\item 
If $k_1^2=\beta^2/4-\gamma\not\equiv 0$ is a constant and $A=0$ then
$w(z)=
\pm k_1 z+c_1-\frac 12\int\beta \,{\rm d}z$.
\item
If $\beta^2/4-\gamma\equiv 0$ then
$\displaystyle
w={\rm e}^{-Az/2}\left\{
c_1-\int\frac\beta 2{\rm e}^{Az/2}{\rm d}z
\right\}$.
\end{enumerate}
%\item
%If $\alpha$ and $\gamma$ are constants and $\beta\equiv 0$, then
%$w(z)=\frac\alpha2(z-c_1)^2+???$.
%\item
\end{enumerate}
%}
\end{theorem}

%\section{Some lemmas}

\vskip 2mm

We have used $k_1$ and $k_2$ to denote constants that appear in constraints on the coefficient functions.  The constants $c_1$ and $c_2$ are parameters in families of solutions of Eq.~\eqref{he}, i.e., they are integration constants.

The case in which $\alpha$, $\beta$ and $\gamma$ are constants was
solved in \cite{chiangh:03}.  In this case any non-constant solution
is  admissible.  Since it is trivial to find the constant solutions,
all meromrophic solutions were found.  

In \cite{hayman:96}, Hayman
conjectured that all entire solutions of
\begin{equation}
\label{original} ff''-f'^2=\kappa_0+\kappa_1f+\kappa_2f'+\kappa_3f''
\end{equation}
have finite order, where $\kappa_0,\ldots,\kappa_3$ are rational functions of
$z$. If we let $w=f-\kappa_3$, then $w$ solves Eq.~\eqref{he} with
$\alpha=\kappa_1-\kappa_3''$, $\beta= \kappa_2+\kappa_3'$ and
$\gamma=\kappa_0+\kappa_1\kappa_3+\kappa_2\kappa_3'+(\kappa_3')^2$. This provided the initial
motivation for studying the meromorphic solutions of
Eq.~\eqref{he}.  However, the problem of the explicit determination of all meromorphic solutions
soon became the main problem of interest.  Nevertheless, Hayman's question is answered by the following elementary corollary of Theorem \ref{mainthm}.

\begin{corollary}
If $\alpha$, $\beta$ and $\gamma$ are rational functions then any transcendental meromorphic solution $w$ of Eq.~(\ref{he}) is of order one, exponential type. 
\end{corollary}

\noindent
This corollary follows immediately on noting that any {\em meromorphic} function that can be expressed as an integral of the form
$\int\beta {\rm e}^{Az}$, for some constant $A$, is itself of the form $B(z){\rm e}^{Az}$, for some rational function $B$.  This can be seen by 
decomposing $\beta$ into partial fractions, using integration by parts, and noting that the coefficients of terms of the form $\int (z-c)^{-1}{\rm e}^{Az}$, where $c$ is constant, must vanish.  In \cite{barsegianll:08},  Barsegian, Laine and L\^e obtained some estimates for the number of poles of meromorphic solutions of Eq.~(\ref{he}) in the case in which $\alpha$, $\beta$ and $\gamma$ are polynomials.

\vskip 2mm

In some sense, Eq.~\eqref{original} is
the simplest differential equation which is neither covered by the
results of Steinmetz (\cite{steinmetz}, \cite[Theorem 12.2]{Laine})
nor Hayman \cite[Theorem C]{hayman:96}. Both these results
 generalise  the classical Gol'dberg Theorem
\cite{Goldberg} that all meromorphic solutions of the first-order
ODE $\Omega(z,f,f')=0$, where $\Omega$ is polynomial in all its
arguments, are of finite order.

Eq.~\eqref{he} is singular when $w=0$.  Suppose that $w$ has a zero at $z=z_0$, which is neither a zero nor a pole of the coefficients, and substitute the expansion
$$
w(z)=\sum_{n=0}^\infty a_n(z-z_0)^{n+p}
$$
in Eq.~\eqref{he}, where $a_0\ne 0$ and $p$ is a positive integer.  If $\gamma\not\equiv 0$ then $p=1$ and there are (generally) two possible values for $a_0$ given by
$a_0^2+\beta(z_0)a_0+\gamma(z_0)=0$.  For each choice of $a_0$ we have a recurrence relation of the form
\begin{equation}
\label{recurrence}
(n+1)(n-r)a_0a_n=P_n(a_0,\ldots,a_{n-1}),
\end{equation}
where for each $n$, $P_n$ is a polynomial in its arguments.  For Eq.~\eqref{he}, $r$ depends on $\alpha$, $\beta$, $\gamma$ and $a_0$.

If $r$ is not a positive integer then all of the coefficients $a_n$ are determined by the choice of $a_0$. % If $r$ is not a positive integer for each $a_0$ then 
In this case there are at most two solutions with a zero at $z_0$.  This is the so-called {\em finiteness property} that has been used by several authors to characterise meromorphic solutions of equations \cite{hille:78,eremenko:82,eremenko:06,eremenkoln:09,contenw:12}.  It is particularly effective for constant coefficient equations as it can be used to deduce periodicity of solutions.

If $r$ is a positive integer then only $a_1,\ldots,a_{r-1}$ are determined by $a_0$.  Eq.~(\ref{recurrence}) shows that there is a necessary
({\em resonance}) condition, $P(a_0,\ldots,a_{r-1})=0$, which must be satisfied.  Subject to this constraint, all remaining coefficients, $a_{r+1},\ldots$ are determined by $a_r$ and $a_0$.  This is very useful for identifying equations that admit meromorphic solutions (see, e.g., \cite{wittich:53}).  One of the main difficulties with Eq.~\eqref{he} is that the location of the resonance depends on the coefficients: $r=(\beta(z_0)/a_0)+2$.  So even in the constant coefficient case considered in \cite{chiangh:03}, we can choose $\beta$ and $\gamma$ so that there is a 
positive integer
resonance at an arbitrary high coefficient in the expansion for $w$,
implying that high-order derivatives of $w$ at a zero of $w$ are not determined by the equation and leading-order term (c.f. \cite{liao:11}).

In the present paper we bypass issues related to resonance by using at most the first two terms in the series expansion for $w$ at zeros to construct a small (in the sense of Nevanlinna theory) function of $w$ and $w'$, the coefficient functions $\alpha$, $\beta$, $\gamma$ and their derivatives.  In this way we construct first-order equations that we can solve for $w$.

\section{Proof of Theorem \ref{mainthm}}
\label{s-proof} If $\alpha(z)\equiv\beta(z)\equiv \gamma(z)\equiv 0$
then Eq.~\eqref{he} becomes $(w'/w)'=0$, which has the general
solution $w(z)=c_2{\rm e}^{c_1z}$. 
This is a special case of part \ref{concC} of Theorem  \ref{mainthm}.
From now on we take at least one
of $\alpha$, $\beta$, $\gamma$ to be nonzero.

For any meromorphic function $f$, we define the set $\Phi_f$ as
follows.  If $f\equiv 0$ then $\Phi_f=\emptyset$.  If $f\not\equiv
0$ then $\Phi_f$ is the set of all zeros and poles of $f$. Let
$\Phi=\Phi_\alpha\cup\Phi_\beta\cup\Phi_\gamma$.
Let
$w\not\equiv 0$ be a meromorphic solution of Eq.~\eqref{he} and let
$z_0\in\Omega:=\mathbb{C}\setminus\Phi$ be either a zero or a pole
of $w$.  Then $w$ has a Laurent series expansion of the form
$$
w(z)=a_0\zeta^p+a_1\zeta^{p+1}+O(\zeta^{p+2}),
$$
where $\zeta=z-z_0$, $a_0\ne 0$ and $p\in\mathbb{Z}\setminus\{0\}$.
Eq.~\eqref{he} then becomes
\begin{equation}
\label{balance}
-pa_0^2\zeta^{2p-2}+\cdots=\alpha(z_0)(a_0\zeta^p+\cdots)+\beta(z_0)(a_0p\zeta^{p-1}+\cdots)+(\gamma(z_0)+\cdots).
\end{equation}
It follows that if $\beta\equiv \gamma\equiv 0$, then $p=2$. Otherwise
$p=1$.  In particular, $w$ is analytic on $\Omega$.\

Throughout this proof we will use the standard notation from
Nevanlinna theory (see, e.g., Hayman \cite{hayman:64} or Laine
\cite{Laine}). In particular, for any meromorphic function $f$, we
denote the (integrated) counting function with multiplicities by
$N(r,f)$ and without multiplicities by $\overline N(r,f)$.
Furthermore we will denote by $N_\Phi(r,f)$ and
$\overline{N}_\Phi(r,f)$ the counting functions (with and without
multiplicities respectively) where we only count the poles of $f$ in
the set $\Phi$. In particular it follows that if $w$ is a
meromorphic solution of Eq.~\eqref{he} then $N(r,w)=N_\Phi(r,w)$.
Now for any meromorphic function $f$, $\overline N_\Phi(r,f)\le
\overline N(r,\alpha)+ \overline N(r,1/\alpha)+  \overline
N(r,\beta)+ \overline N(r,1/\beta)+ \overline N(r,\gamma)+ \overline
N(r,1/\gamma)=S(r,w)$, where if $\alpha\equiv 0$ we take
$\overline N(r,1/\alpha)=0$, etc. So $\overline N(r,w)=\overline
N_\Phi(r,w)=S(r,w)$.

When the coefficient functions $\alpha$, $\beta$ and $\gamma$ are rational functions then $\Phi$ is a finite set and $N_\Phi(r,w)=S(r,w)$.  However, for transcendental coefficients this does not follow immediately.

\medskip

\noindent{\bf Case 1:} $\alpha\not\equiv 0$, $\beta\equiv\gamma\equiv 0$.\\
Substituting $w(z)=a_0\zeta^2+a_1\zeta^3+O(\zeta^4)$ in
Eq.~\eqref{he} shows that about any $z_0\in\Omega$ such that
$w(z_0)=0$, we have
$$
w(z)=-\frac{\alpha(z_0)}{2}\zeta^2-\frac{\alpha'(z_0)}{2}\zeta^3+O(\zeta^4).
$$
Together with the fact that $w$ is analytic on $\Omega$, it follows that
\begin{equation}
\label{f1}
f(z):=\left(
\frac{w'}{w}-\frac{\alpha'}{\alpha}
\right)^2+2\frac\alpha w
\end{equation}
is also analytic on $\Omega$.  Using Eq.~\eqref{he} with
$\beta\equiv\gamma\equiv 0$, we see that
\begin{equation}
\label{fa-subs}
f(z)=\left(\frac{w'}{w}-\frac{\alpha'}{\alpha}
\right)^2+2\left(\frac{w'}{w}\right)'.
\end{equation}
Hence
\begin{eqnarray*}
N(r,f)&=&N_\Phi(r,f)\le 2N_\Phi\left(r,\frac{w'}{w}\right)+2N_\Phi\left(r,\frac{\alpha'}{\alpha}\right)+N_\Phi\left(r,\left(\frac{w'}{w}\right)'\right)
\\
&=&
4\left\{\bar N_\Phi(r,w)+\bar N_\Phi\left(r,\frac1w\right)\right\}+2N_\Phi\left(r,\frac{\alpha'}{\alpha}\right)=S(r,w).
\end{eqnarray*}
Furthermore, applying the Lemma on the Logarithmic Derivative to
Eq.~\eqref{fa-subs} gives $m(r,f)=S(r,w)$.  So $T(r,f)=S(r,w)$.

Differentiating Eq.~\eqref{f1} and using Eq.~\eqref{he} to eliminate $w''$ gives
\begin{equation}
\label{f1p}
%\begin{split}
f'
%&=2\Big(\frac{w'}{w}-\frac{\alpha'}{\alpha}\Big)\Big[\Big(\frac{w'}{w}\Big)'-\Big(\frac{\alpha'}{\alpha}\Big)'\Big]+2\frac{\alpha'}{w}-2\alpha\frac{w'}{w^2}\\
%&=2\Big(\frac{w'}{w}-\frac{\alpha'}{\alpha}\Big)\Big[\frac{\alpha}{w}-\Big(\frac{\alpha'}{\alpha}\Big)'\Big]+2\frac{\alpha'}{w}-2\alpha\frac{w'}{w^2}\\
%&
=2\Big(\frac{\alpha'}{\alpha}\Big)'\Big(\frac{\alpha'}{\alpha}-\frac{w'}{w}\Big).
%\end{split}
\end{equation}
When $(\frac{\alpha'}{\alpha})'\not\equiv0$, we obtain
\begin{equation*}
\frac{w'}{w}=\frac{\alpha'}{\alpha}-\frac{f'}{2}\Big[\Big(\frac{\alpha'}{\alpha}\Big)'\Big]^{-1}.
\end{equation*}
Substituting this into Eq.~\eqref{f1} gives
\begin{equation}\label{f,w}
f=\frac{f'^2}{4}\Big[\Big(\frac{\alpha'}{\alpha}\Big)'\Big]^{-2}+2\frac{\alpha}{w}.
\end{equation}
Applying Nevanlinna's First Fundamental Theorem to Eq.~\eqref{f,w}, we
obtain $T(r,w)=S(r,w)$, a contradiction. Therefore
$(\alpha'/\alpha)'\equiv0$, so from Eq.~\eqref{f1p}, $f'\equiv 0$.
Thus,  \begin{equation*}\alpha(z)\equiv k_1e^{k_2z}\quad
\text{and}\quad f(z)\equiv c_1^2,\end{equation*} where $k_1\not=0$,
$k_2$ and $c_1$ are constants.

In terms of $u=w/\alpha$, Eq.~\eqref{f1} becomes
$$
u'^2=c_1^2u^2-2u.
$$
If $c_1=0$, this gives $u=-\frac 12(z+c_2)^2$. When $k_2\not=0$, we
arrive at the contradiction $T(r,w)=T(r,\alpha)+S(r,\alpha)=S(r,w)$.  Thus $k_2=0$, so $\alpha=k_1$ and
$w=-\frac{k_1}{2}(z+c_2)^2$. For the case $c_1\not=0$,
$u=c_1^{-2}\{\cosh(c_1z+c_2)+1\}$ where $c_2$ is a constant. This gives part \ref{concA} of Theorem \ref{mainthm}.

\medskip

\noindent{\bf Case 2:} $\beta\not\equiv0$, $\gamma\equiv 0$.\\
Recall that in this case $w$ has only simple zeros in $\Omega$.
Substituting $w(z)=a_0\zeta+a_1\zeta^2+O(\zeta^3)$ in Eq.~\eqref{he}
yields, at the leading order $a_0=-\beta(z_0)$ and at the
next-to-leading order we find the constraint
$\alpha(z_0)+\beta'(z_0)=0 $.

\medskip

\noindent{\bf Case 2a:} $\alpha+\beta'\not\equiv 0$, $\gamma\equiv 0$.\\
Let $f=w'/w$.  If $z_0$ is a pole of $w$ then $z_0\in\Phi$.  If
$z_0$ is a zero of $w$, then either $z_0\in\Phi$ or
$\alpha(z_0)+\beta'(z_0)=0$.  Hence
$$
N(r,f)=N\left(r,\frac{w'}w\right)\le\overline N_\Phi(r,w)+ \overline N_\Phi\left(r,\frac1w\right)+\overline N\left(r,\frac1{\alpha+\beta'}\right)=S(r,w).
$$
 It then follows form the Lemma on the Logarithmic Derivative that $T(r,f)=T(r,w'/w)=S(r,w)$.  Substituting $w'=fw$ and $w''=(f'+f^2)w$ in Eq.~\eqref{he} with $\gamma\equiv 0$ yields
 $$
 f'w=\alpha+f\beta.
 $$
Since $T(r,f')=S(r,w)$ and $T(r,\alpha+f\beta)=S(r,w)$, we must have $f'\equiv\alpha+f\beta\equiv 0$, thus $f\equiv k_1$ is a constant.
Hence there exists a constants $k_1$ such that
$w(z)=c_1{\rm e}^{k_1 z}$ and $\alpha(z)=-k_1\beta(z)$, giving part \ref{concB} of the theorem.

\medskip

\noindent{\bf Case 2b:} $\alpha+\beta'\equiv 0$, $\gamma\equiv 0$.\\
Eq.~\eqref{he} takes the form
$\left(({w'+\beta})/w\right)'=0$,
%hence
%\begin{equation}
%\label{first-2b}
%w'=c_1w-\beta,
%\end{equation}
%for some constant $c_1$. % It follows that $w''=c_1^2w-(\beta'+c_1\beta)$ and so any solution of equation Eq.\eqref{first-2b} also solves Eq.\eqref{he} with $\alpha+\beta'\equiv\gamma\equiv 0$. By the fundamental fact of first order linear ODE, we have 
which has the general solution $w=e^{c_1z}\{c_2-\int \beta(z)e^{-c_1z}dz\}$ where $c_2$ is a constant.  This gives part \ref{concC} of the theorem.

\medskip

\noindent{\bf Case 3:} $\gamma\not\equiv 0$.\\
Recall that in this case $w$ is analytic in $\Omega$ and any zero $z_0$ of $w$ in $\Omega$ is simple.  On substituting $w(z)=a_0\zeta+a_1\zeta^2+O(\zeta^3)$ in Eq.~\eqref{he} we find that
$$
a_0^2+\beta(z_0)a_0+\gamma(z_0)=0\quad\mbox{and}\quad
%$$
%and
%$$
a_1=\frac{1}{2\gamma(z_0)}\left\{
\gamma'(z_0)-\beta(z_0)(\alpha(z_0)+\beta'(z_0))
\right\}a_0-\frac12(\alpha(z_0)+\beta'(z_0)).
$$
Let
$$
f(z)=\frac{(w')^2+\beta w'+\gamma}{w^2}.
$$
If $f$ has a pole at $z_0\in\Omega$ then $w(z_0)=0$, 
From Eq.~(\ref{he}), $f(z)=(w''-\alpha)/w$,
so in a neighbourhood of $z_0$,
$$
f(z)=\frac{2a_1-\alpha(z_0)}{a_0\zeta}+O(1)
=
\left\{
\frac{\gamma'(z_0)-\beta(z_0)[\alpha(z_0)+\beta'(z_0)]}{\gamma(z_0)\zeta}-\frac{2\alpha(z_0)+\beta'(z_0)}{a_0\zeta}
\right\}
+O(1).
$$
Therefore
\begin{equation}
\label{gdef}
g(z)=\frac{(w')^2+\beta w'+\gamma}{w^2}+A\frac{w'}w+\frac {2\alpha+\beta'}w,\quad A=\frac{\beta(\alpha+\beta')-\gamma'}{\gamma},
\end{equation}
is analytic on $\Omega$.

Rewriting Eq.~\eqref{he} as
\begin{equation}
\label{rewrite}
\frac 1{w^2}=\frac{1}{\gamma}\left\{
\left(\frac{w'}{w}\right)'-\frac{1}{w}\left(\alpha+\beta\frac{w'}w\right)
\right\},
\end{equation}
we see that
$$
2N_\Phi\left(r,\frac1{w}\right)
\le
N_\Phi\left(r,\frac1\gamma\right)+N_\Phi\left(r,\left(\frac{w'}w\right)'\right)+N_\Phi\left(r,\frac 1w\right)+N_\Phi\left(r,\alpha\right)+N_\Phi\left(r,\beta\right)+N_\Phi\left(r,\frac{w'}w\right).
$$
Hence
$$
N_\Phi\left(r,\frac1{w}\right)\le 3\left\{ \overline N_\Phi\left(r,{w}\right)+\overline N_\Phi\left(r,\frac1{w}\right)\right\}+S(r,w)=S(r,w).
$$
So from Eq.~\eqref{gdef}, we have
$$N(r,g)=N_\Phi(r,g)\leq 2\overline{N}_\Phi(r,w)+2N_\Phi\left(r,\frac{1}{w}\right)+S(r,w)=S(r,w).$$
Taking the proximity function of both sides of Eq.~\eqref{rewrite}, we obtain
\begin{equation*}
\begin{split}
2m\left(r,\frac1w\right)
&\le
m\left(r,\frac1\gamma\right)+m\left(r,\left(\frac{w'}w\right)'\right)+m\left(r,\frac1w\right)+m\left(r,\alpha\right)+m\left(r,\beta\right)+m\left(r,\frac{w'}w\right)
\\
&= m\left(r,\frac1w\right)+S(r,w).
\end{split}
\end{equation*}
Hence $m(r,g)=S(r,w)$.  So $T(r,g)=m(r,g)+N(r,g)=S(r,w)$.

Differentiating $w^2\times\,$Eq.~\eqref{gdef} and using Eq.~\eqref{he} to eliminate $w''$ and Eq. \eqref{gdef} to eliminate $(w')^3$ and then $(w')^2$, we have
%\begin{equation}\label{gprime}
%w^3g'=\beta(w')^2+(A'w^2+\beta A w+\beta^2)w'+\beta\gamma+\beta(2\alpha+\beta ')w+(B-\beta g)w^2,
%\end{equation}
%(Eq.~\ref{gprime})$-\beta w^2\times$(Eq.~\ref{gdef}) gives
\begin{equation}
\label{lastwp}
A'w'=g'w-B,
\end{equation}
where $B=\beta g+\alpha A+2\alpha'+\beta''$.

\medskip

\noindent{\bf Case 3a:}  $A'\not\equiv 0$.\\
Using Eq.~\eqref{lastwp} to eliminate $w'$ from Eq.~\eqref{gdef} gives
$$
(g'^2-gA'^2+g'AA')w^2+([2\alpha+\beta']A'^2+\beta g' A'-AA'B-2g'B)w+(B^2-\beta A'B+\gamma A'^2)=0.
$$
Since the coefficients of the different powers of $w$ are all $S(r,w)$, we must have that each coefficient vanishes identically.  In particular, the coefficient of $w^2$ gives
%\begin{equation}
%\label{gA}
%\begin{split}
$g'^2-gA'^2+g'AA' =0$.
%,\\
%(2\alpha+\beta')A'^2+\beta g' A'-AA'B-2g'B &=0,\\
%B^2-\beta A'B+\gamma A'^2 &=0.
%\end{split}
%\end{equation}
%which has the general solution
It follows that $G=g+(A^2/4)$ satisfies $(G')^2=(A')^2G$.  Hence either $G=0$ (i.e., $g=-A^2/4$) or $G=(A/2+k_1)^2$ (i.e., $g=k_1A+k_1^2$), where $k_1$ is a constant.

\medskip

\noindent{\bf Case 3a(i):}
$g=k_1A+k_1^2$.

\noindent
Eq.~\eqref{lastwp} now has the form
\begin{equation}
\label{wpc}
w'=k_1w+h,
\end{equation}
where $h=-B/A'$.  Hence $w''=k_1^2w+(h'+k_1h)$ and we see that any solution of Eq.~\eqref{wpc} solves Eq.~\eqref{he} if and only if
$$(h'-k_1h-\alpha-k_1\beta)w=h^2+\beta h+\gamma,$$
so $h^2+\beta h+\gamma=0$
and
$h'-k_1h=\alpha+k_1\beta$.  This corresponds to part \ref{concD} of the theorem.
\medskip

\noindent{\bf Case 3a(ii):} $g=-A^2/4$.

\noindent
Eq.~\eqref{lastwp} becomes $w'=-(A/2)w+h$, where $h=-B/A'$.  Hence $w''=[(A^2/4)-(A'/2)]w+h'-hA/2$.
Using these expressions to eliminate the first and second derivatives in Eq.~\eqref{he} leads to
$$-\frac{A'}{2}w^2+\left(\frac{Ah}{2}+h'-\alpha+\frac{\beta A}{2}\right)w=h^2+\beta h+\gamma$$
 with coefficients that are $S(r,w)$. By the Valiron-Mokhon'ko Theorem \cite[Theorem 1.13]{yangyi}, we have $2T(r,w)=S(r,w)$, which is impossible.

\medskip

\noindent{\bf Case 3b:}   $A'\equiv 0$, i.e. $A$ is a constant.\\
It follows from Eq.~\eqref{lastwp} that $g$ is also a constant and $B=0$.  Eq.~\eqref{gdef} can be rewritten as
\begin{equation}
\label{Ap0}
(w'+\frac 12[Aw+\beta])^2
%=gw^2-(2\alpha+\beta')w-\gamma+\frac14\left(Aw+\beta\right)^2
=
\left(g+\frac{A^2}4\right)w^2+\left(
\frac\beta 2A-\beta'-2\alpha
\right)w
+
\left(\frac{\beta^2}4-\gamma\right).
\end{equation}
Let
$h(z)=\left(
\frac\beta 2A-\beta'-2\alpha
\right){\rm e}^{Az/2}$.
Then
\begin{equation}
\label{Discp}
\left(\left[\frac{\beta^2}4-\gamma\right]{\rm e}^{Az}\right)'=\frac\beta 2{\rm e}^{Az/2}h
\end{equation}
and the condition $B=0$ is equivalent to
\begin{equation}
\label{hp}
h'=\left(g+\frac{A^2}4\right)\beta{\rm e}^{Az/2}.
\end{equation}
Clearly if $g=-A^2/4$ then $h$ is constant.

\medskip

\noindent{\bf Case 3b{(i)}:}  $g+\frac{A^2}4\ne 0$.
\\
So $k_1^2=g+\frac{A^2}4$ is a non-zero constant.
It follows from Eqs.\eqref{Discp} and \eqref{hp} that
$$\left(\left[\frac{\beta^2}4-\gamma\right]{\rm e}^{Az}\right)'=\frac{1}{2k_1^2}hh'.
$$
Integration shows that
\begin{equation}
\label{k1,k2}
k_2^2=\frac1{k_1^2}
\left\{
\frac{h^2}{4k_1^2}+\left(\gamma-\frac{\beta^2}4\right){\rm e}^{Az}
\right\}
=
\frac1{k_1^2}
\left\{
\frac{1}{4k_1^2}\left(
\frac\beta 2A-\beta'-2\alpha
\right)^2+\left(\gamma-\frac{\beta^2}4\right)
\right\}{\rm e}^{Az}
\end{equation}
is a constant. Let
$$
u=w{\rm e}^{Az/2}+\frac{h}{2k_1^2}.
$$
Then Eq.~\eqref{Ap0} becomes \begin{equation}\label{simple}(u')^2=k_1^2(u^2-k_2^2).\end{equation} 
When $k_2\ne 0$ we have 
$$w=\left(\pm k_2\cosh(k_1z+c_1)-\frac{h}{2k_1^2}\right){\rm e}^{-Az/2},$$
where $c_1$ is a constant.  Therefore $T(r,w)\le K_1r+S(r,w)$ for some $K_1>0$.
When $A\ne 0$,  Eq.~\eqref{k1,k2} shows that
$r\le K_2T(r,{\rm e}^{Az})=S(r,w)$, which gives the contradiction $T(r,w)=S(r,w)$.  Hence $A=0$ if $k_2\ne 0$.
This is part \ref{concE}(a) of the theorem.  Part \ref{concE}(b) corresponds to the case in which
$k_2=0$, where
$$
w=e^{-Az/2}\left(c_1e^{\pm k_1z}-\frac{h}{2k_1^2}\right)=c_1{\rm e}^{\left(-\frac A2\pm k_1\right)z}-\frac 1{2k_1^2}\left(\frac\beta 2A-\beta'-2\alpha\right).
$$

\medskip

\noindent{\bf Case 3b{(ii)}:}  $g=-A^2/4$, $h\ne0$.
\\
Let $\lambda=\int\frac\beta 2 {\rm e}^{Az/2}{\rm d}z$.  It follows from Eq.~\eqref{Discp} and Eq.~\eqref{hp} that $h$ and
$$
C:=\frac 1h\left(\frac{\beta^2}4-\gamma\right){\rm e}^{Az}-\lambda
$$
are constants.  Let $u=w{\rm e}^{Az/2}+\lambda$.  Then Eq.~\eqref{Ap0} becomes
$(u')^2=h(u+C)$, which has the general solution
$u=\frac h4(z+c_1)^2-C$.
Hence
\begin{equation}
\label{quad}
w=\frac h4{\rm e}^{-Az/2}(z+c_1)^2-\frac 1h\left(\frac{\beta^2}4-\gamma\right){\rm e}^{Az/2}
=\frac14\left(\frac\beta 2A-\beta'-2\alpha\right)(z+c_1)^2-\frac{\frac{\beta^2}4-\gamma}{\frac\beta 2A-\beta'-2\alpha}.
\end{equation}
So $T(r,w)=O(r)+S(r,w)$.  Recall that $h$ is a nonzero constant.
Now if $A\ne 0$, we have $r\le K_1 T(r,{\rm e}^{Az/2})=K_1T(r,\frac\beta 2A-\beta'-2\alpha)+O(1)=S(r,w)$, a contradiction.
Therefore $A=0$.  Now Eq.~\eqref{quad} with $A=0$ shows that $T(r,w)=2\log r+S(r,w)$.  Hence $w$ is admissible if and only if the coefficients
$\alpha$, $\beta$ and $\gamma$ are constants.
This gives part \ref{concE}(c).

\medskip

\noindent{\bf Case 3b{(iii)}:}  $g=-A^2/4$, $h=0$.
\\
It follows from Eq.~\eqref{Ap0} and Eq.~\eqref{Discp} that
$w'+\frac 12[Aw+\beta]=k_1{\rm e}^{-Az/2}$,
where
$k_1^2=\left(\frac{\beta^2}4-\gamma\right){\rm e}^{Az}$
is a constant. Hence $\left(w{\rm e}^{Az/2}\right)'=k_1-\frac 12\beta{\rm e}^{Az/2}$, giving
$$
w={\rm e}^{-Az/2}\left\{
k_1 z+c_1-\int\frac\beta 2{\rm e}^{Az/2}{\rm d}z
\right\}.
$$
To study the admissibility of this solution, we will use the following theorem from Hayman and  Miles \cite{haymanm:89}.

\begin{theorem}
Let $f(z)$ be a transcendental meromorphic function and $K>1$ be a real number. Then
there exists a set $M(K)$ of upper logarithmic density at most
$
d(K)=1-(2{\rm e}^{K-1}-1)^{-1}>0
$
such that for every positive integer $k$, we have
$$
\limsup_{r\to\infty,\ r\not\in M(K)}\frac{T(r,f)}{T(r,f^{(k)})}\le 3{\rm e}K.
$$
\end{theorem}

Furthermore, note that if $f$ is any non-constant rational function other than a degree one polynomial, then $T(r,f)\le K T(r,f')$ for some $K>0$.
Therefore if $w{\rm e}^{Az/2}$ is not a constant or a degree one polynomial and $k_1\ne 0$, it follows that there is a sequence of values of $r\to\infty$ such that for some $K_1>0$,
\begin{equation*}
\begin{split}
T(r,w) &\le T\left(r,w{\rm e}^{Az/2}\right) + T\left(r,{\rm e}^{-Az/2}\right)\\
&\le K_1 T\left(r,\left(w{\rm e}^{Az/2}\right)'\right) + T\left(r,{\rm e}^{-Az/2}\right)\\
&= K_1 T\left(r,\frac 12\beta{\rm e}^{Az/2}-k_1\right) + T\left(r,{\rm e}^{-Az/2}\right)\\
&= K_1 T\left(r,\frac 12k_1\beta\left(\frac{\beta^2}4-\gamma\right)^{-1/2}-k_1\right) + T\left(r,k_1^{-1}\left(\frac{\beta^2}4-\gamma\right)^{1/2}\right)\\
&=o(T(r,w)),
\end{split}
\end{equation*}
which is a contradiction.  If $w{\rm e}^{Az/2}$ is at most a degree one polynomial, then $\beta=k_2{\rm e}^{-Az/2}$ and $w$ is only admissible if $A=0$.
Now $w$ is a polynomial of degree no more than one, so
$\alpha$, $\beta$ and $\gamma$ are constants. It follows from $h=0$ that $\alpha=0$.  At the same time, $A=0$
and $\alpha=0$ implies that $g=0$, so we have 
$w'^2+\beta w'+\gamma=0$. This corresponds to
part  \ref{concE}(d) of the theorem.
Otherwise we have $k_1=0$, i.e. $\gamma=\frac{\beta^2}4$, which corresponds to part \ref{concE}(e).

\section{Discussion}
The proof provided in Section \ref{s-proof} would have been significantly shorter had we restricted ourselves to the rational coefficient case.  In the first instance, the fact that $N(r,1/w)=\overline N(r,1/w)+S(r,w)$ would have followed immediately from Eq.~\eqref{balance}.  Also, many of the subcases considered in the proof could be eliminated or simplified because they require that in general certain rational functions of the coefficient functions be an exponential in $z$.

When we allowed some of the coefficients to be transcendental, we generated many exact solutions only to discard them later because these solutions grow at the same rate as the coefficients.  From the point of view of using the existence of meromorphic solutions as a detector of exactly solvable cases, this suggests that 
perhaps a weaker
 notion of ``admissibility'' would be more fruitful.  These are all perfectly good solutions and it is undesirable  merely  to discard them or even to search for more efficient methods to avoid considering them in the first place.  It seems wasteful not to modify the problem so that such solutions will appear in the final classification.  We hope to explore this problem  in future work.

%Why we need a better notion than ``admissibility''.

\vskip 5mm

\noindent{\bf \large Acknowledgements}\vskip 4mm
The first author was partially supported by the Engineering and Physical Sciences Research Council grant numbers EP/I013334/1 and EP/K041266/1.
The second author was partially supported by an Erasmus Mundus scholarship through the EMEA scholarship programme and the National Science Foundation of China No.~11001057.

\end{document}